\documentclass[12pt]{article}
\usepackage{amsmath}
\usepackage{amssymb}
\usepackage{amsfonts}
\usepackage{amsthm}
\usepackage{enumerate}

\def\init{\setcounter{equation}{0}}

\newtheorem{theorem}{Theorem}[section]
\newtheorem{proposition}[theorem]{Proposition}
\newtheorem{lemma}[theorem]{Lemma}

\newtheorem{corollary}[theorem]{Corollary}

\def\Res{{\operatorname{Res}}}
\def\m{{\operatorname{m}}}

\def\d{{\operatorname{d}}}

\def\Re{\hbox{Re}\,}

\def\ii{{\mathrm{i}}}
\def\e{{\mathrm{e}}}

\numberwithin{equation}{section}
\newenvironment{acknowledgements}{\noindent{\bf Acknowledgements}\bigskip}{}

\begin{document}

\title{New $_5F_4$ hypergeometric transformations, three-variable Mahler
measures, and formulas for $1/\pi$}
\author{Mathew D. Rogers \\
        \small{\textit{Department of Mathematics, University of British
        Columbia}}\\
        \small{\textit{Vancouver, BC, V6T-1Z2, Canada}}\\
        \small{\textbf{email:} matrogers@math.ubc.ca}}

\maketitle

\abstract{New relations are established between families of
three-variable Mahler measures.  Those identities are then expressed
as transformations for the $_5F_4$ hypergeometric function.  We use
these results to obtain two explicit $_5F_4$ evaluations, and
several new formulas for $1/\pi$.}

\medskip


\section{Introduction}
\label{intro} \init
In this paper we will study the consequences of some recent results
of Bertin.  Recall that Bertin proved $q$-series expansions for a
pair of three-variable Mahler measures in \cite{bert2}.  As usual
the Mahler measure of an $n$-variable polynomial,
$P(z_1,\dots,z_n)$, is defined by
\begin{equation*}
\m\left(P(z_1,\dots,z_n)\right)=\int_{0}^{1}\dots\int_{0}^{1}\log\left|P\left(\e^{2\pi\ii\theta_1},\dots,\e^{2\pi\ii\theta_n}\right)\right|\d\theta_1\dots\d\theta_n.
\end{equation*}
We will define $g_1(u)$ and $g_2(u)$ in terms of the following
three-variable Mahler measures
\begin{align}
g_1(u):=\m&\left(u+x+\frac{1}{x}+y+\frac{1}{y}+z+\frac{1}{z}\right),\label{bertin
g1 definition}\\
\begin{split}
g_2(u):=\m&\left(-u+4+\left(x+x^{-1}\right)\left(y+y^{-1}\right)\right.\\
&\left.\quad+\left(x+x^{-1}\right)\left(z+z^{-1}\right)+\left(y+y^{-1}\right)\left(z+z^{-1}\right)\right).\label{bertin
g2 definition}
\end{split}
\end{align}
We can recover Bertin's original notation by observing that
$g_1(u)=\m\left(P_u\right)$, and after substituting
$(xz,y/z,z/x)\rightarrow(x,y,z)$ in Eq. \eqref{bertin g2 definition}
we see that $g_2(u+4)=\m\left(Q_u\right)$ \cite{bert2}.

In Section \ref{identities between Mahler measures} we will show how
to establish a large number of interesting relations between
$g_1(u)$, $g_2(u)$, and three more
three-variable Mahler measures. 
For example, for $|u|$ sufficiently large Eq. \eqref{g1 in terms of
f4 explicit} is equivalent to
\begin{equation}\label{intro g1 and f4 identity}
\begin{split}
g_1\left(3\left(u^2+u^{-2}\right)\right)=&\frac{1}{5}\m\left(x^4+y^4+z^4+1+\sqrt{3}\frac{\left(3+u^4\right)}{u^{3}}x
y
z\right)\\
&+\frac{3}{5}\m\left(x^4+y^4+z^4+1+\sqrt{3}\frac{\left(3+u^{-4}\right)}{u^{-3}}x
y z\right).
\end{split}
\end{equation}
Rodriguez-Villegas briefly mentioned the Mahler measure
$\m\left(x^4+y^4+z^4+1+u x y z\right)$ on the last page of
\cite{RV}.

We will also show that identities like Eq. \eqref{intro g1 and f4
identity} are equivalent to transformations for the $_5F_4$
hypergeometric function.  Recall that the generalized hypergeometric
function is defined by
\begin{equation*}
{_pF_q}\left(\substack{a_1,a_2,\dots a_p\\b_1,b_2,\dots
b_q};x\right)=\sum_{n=0}^{\infty}\frac{(a_1)_n\dots(a_p)_n}{(b_1)_n\dots(b_q)_n}\frac{x^n}{n!},
\end{equation*}
where $(y)_n=\Gamma(y+n)/\Gamma(y)$.  We have restated Eq.
\eqref{intro g1 and f4 identity} as a hypergeometric transformation
in Eq. \eqref{5F4 identity two}. As a special case of Eq.
\eqref{intro g1 and f4 identity} we can also deduce that
\begin{equation*}
{_5F_4}\left(\substack{\frac{5}{4},\frac{3}{2},\frac{7}{4},1,1\\2,2,2,2};1\right)=\frac{256}{3}\log(2)-\frac{5120\sqrt{2}}{3\pi^3}L(f,3),
\end{equation*}
where $
f(q)=q\prod_{n=1}^{\infty}\left(1-q^n\right)^2\left(1-q^{2n}\right)\left(1-q^{4n}\right)\left(1-q^{8n}\right)^2$,
and $L(f,s)$ is the usual $L$-series of $f(q)$.  We will conclude
Section \ref{identities between Mahler measures} with a brief
discussion of some related, but still unproven, evaluations of the
$_4F_3$ and $_3F_2$ hypergeometric functions.

It turns out that $g_1(u)$ and $g_2(u)$ are also closely related to
Watson's triple integrals.  For appropriate values of $u$, Watson
showed that $g_1'(u)$ and $g_2'(u)$ reduce to products of elliptical
integrals (for relevant results see \cite{Wa}, \cite{JZ}, \cite{Jc},
and \cite{Gl}). In Section \ref{pi section} we will use some related
transformations to prove new formulas for $1/\pi$.  For example, we
will show that
\begin{equation*}
\frac{2}{\pi}=\sum_{n=0}^{\infty}(-1)^n\frac{(3n+1)}{32^n}\sum_{k=0}^{n}{2n-2k\choose
n-k}{2k\choose k}{n\choose k}^2.
\end{equation*}
Notice that this formula for $1/\pi$ involves the Domb numbers.
Chan, Chan and Liu obtained a similar formula for $1/\pi$ involving
Domb numbers in \cite{ccl}, we have recovered their result in
\eqref{pi 2}. Zudilin and Yang also discovered some related formulas
for $1/\pi$ in \cite{Zd}. All of the $_3F_2$ transformations that we
will utilize in Section \ref{pi section} follow from differentiating
the $_5F_4$ identities established in Section \ref{identities
between Mahler measures}.

\section{Identities between Mahler measures and transformations for the $_5F_4$ function}
\label{identities between Mahler measures} \init

Bertin proved that both $g_1(u)$ and $g_2(u)$ have convenient
$q$-series expansions when $u$ is parameterized correctly. Before
stating her theorem, we will define some notation.  As usual let
\[(x,q)_{\infty}=(1-x)(1-xq)\left(1-x q^2\right)\dots,\]
and define $G(q)$ by
\begin{equation}\label{G(q) definition}
G(q)=\Re\left[-\log(q)+240\sum_{n=1}^{\infty}n^2\log\left(1-q^n\right)\right].
\end{equation}
Notice that if $q\in(0,1)$ then $G'(q)=-M(q)/q$, where $M(q)$ is the
Eisenstein series of weight $4$ on the full modular group
$\Gamma(1)$ \cite{Be2}.
\begin{theorem}(Bertin)  For $|q|$ sufficiently small
\begin{align}
g_1\left(t_1(q)\right)=&-\frac{1}{60}G(q)+\frac{1}{30}G(q^2)-\frac{1}{20}G(q^3)+\frac{1}{10}G(q^6),\label{bertin
g1 series simplified}\\
g_2\left(t_2(q)\right)=&\frac{1}{120}G(q)-\frac{1}{15}G\left(q^2\right)-\frac{1}{40}G\left(q^3\right)+\frac{1}{5}G\left(q^6\right),\label{bertin
g2 series}
\end{align}
where
\begin{align*}
t_1(q)=&v_1+\frac{1}{v_1}, \text{ and } v_1=q^{1/2}\frac{\left(q;q^2\right)_{\infty}^6}{\left(q^3;q^6\right)_{\infty}^6},\\
t_2(q)=&-\left(v_2-\frac{1}{v_2}\right)^2,\text{ and }
v_2=q^{1/2}\frac{\left(q^2;q^2\right)_{\infty}^6\left(q^3;q^3\right)_{\infty}^2\left(q^{12};q^{12}\right)_{\infty}^4}{\left(q;q\right)_{\infty}^2\left(q^4;q^4\right)_{\infty}^4\left(q^6;q^6\right)_{\infty}^6}.
\end{align*}
\end{theorem}
 In this section we will show that both $g_1(u)$ and $g_2(u)$
reduce to linear combinations of $_5F_4$ hypergeometric functions.
We will accomplish this goal by first expressing each of the
functions
\begin{align}
f_2(u):=&2\m\left(u^{1/2}+\left(x+x^{-1}\right)\left(y+y^{-1}\right)\left(z+z^{-1}\right)\right),\label{f2 definition}\\
f_3(u):=&\m\left(u-\left(x+x^{-1}\right)^2\left(y+y^{-1}\right)^2\left(1+z\right)^3z^{-2}\right),\label{f3 definition}\\
f_4(u):=&4\m\left(x^4+y^4+z^4+1+u^{1/4} x y z\right),\label{f4
definition}
\end{align}
in terms of $G(q)$.  We will then exploit those identities to
establish linear relations between functions in the set
$\{f_2(u),f_3(u),f_4(u),g_1(u),g_2(u)\}$.  This is significant since
$f_2(u)$, $f_3(u)$, and $f_4(u)$ all reduce to $_5F_4$
hypergeometric functions. In particular this implies the non-trivial
fact that both $g_1(u)$ and $g_2(u)$ also reduce to linear
combinations of $_5F_4$ functions. 

\begin{proposition}\label{hypergeometric proposition} The following identities hold for $|u|$ sufficiently large:
\begin{align}
f_2(u)=&\Re\left[\log(u)-\frac{8}{u}{_5F_4}\left(\substack{\frac{3}{2},\frac{3}{2},\frac{3}{2},1,1\\2,2,2,2};\frac{64}{
u}\right)\right],\label{f2 in terms of a 5F4}\\
f_3(u)=&\Re\left[\log(u)-\frac{12}{u}{_5F_4}\left(\substack{\frac{4}{3},\frac{3}{2},\frac{5}{3},1,1\\2,2,2,2};\frac{108}{
u}\right)\right],\label{f3 in terms of a 5F4}\\
f_4(u)=&\Re\left[\log(u)-\frac{24}{u}{_5F_4}\left(\substack{\frac{5}{4},\frac{3}{2},\frac{7}{4},1,1\\2,2,2,2};\frac{256}{
u}\right)\right].\label{f4 in terms of a 5F4}
\end{align}
For $|u|>6$
\begin{align}
g_1(u)=&\Re\left[\log(u)-\sum_{n=1}^{\infty}\frac{\left(1/u\right)^{2n}}{2n}{2n\choose
n}\sum_{k=0}^{n}{2k\choose k}{n\choose k}^2\right],\label{g1 Taylor
series}
\end{align}
and if $|u|>16$
\begin{align}
g_2(u)=&\Re\left[\log(u)-\sum_{n=1}^{\infty}\frac{\left(1/u\right)^{n}}{n}\sum_{k=0}^{n}{{2n-2k\choose
n-k}}{2k\choose k}{n\choose k}^2\right].\label{g2 Taylor series}
\end{align}
\end{proposition}
\begin{proof} We can prove each of these identities using a method due to Rodriguez-Villegas \cite{RV}.
We will illustrate the proof of Eq. \eqref{f2 in terms of a 5F4}
explicitly. Rearranging the Mahler measure shows that
\begin{equation*}
f_2(u)=\Re\left[\log(u)+\int_{0}^{1}\int_{0}^{1}\int_{0}^{1}\log\left(1-\frac{64}{u}\cos^2(2\pi
t_1)\cos^2(2\pi t_2)\cos^2(2\pi t_3)\right)\d t_1\d t_2\d
t_3\right].
\end{equation*}
If $|u|>64$, then $\left|\frac{64}{u}\cos^2(2\pi t_1)\cos^2(2\pi
t_2)\cos^2(2\pi t_3)\right|<1$, hence by the Taylor series for the
logarithm
\begin{equation*}
\begin{split}
f_2(u)=&\Re\left[\log(u)-\sum_{n=1}^{\infty}\frac{(64/u)^{n}}{n}
\left(\int_{0}^{1}\cos^{2n}(2\pi t)\d t\right)^3\right]\\
=&\Re\left[\log(u)-\sum_{n=1}^{\infty}{2n\choose
n}^3\frac{\left(1/u\right)^{n}}{n}\right]\\
=&\Re\left[\log(u)-\frac{8}{u}{_5F_4}\left(\substack{\frac{3}{2},\frac{3}{2},\frac{3}{2},1,1\\2,2,2,2};\frac{64}{
u}\right)\right].
\end{split}
\end{equation*}
Notice that Eq. \eqref{f2 in terms of a 5F4} holds whenever
$u\not\in[-64,64]$, since $f_2(u)$ is harmonic in
$\mathbb{C}\backslash[-64,64]$.$\blacksquare$
\end{proof}

While Proposition \ref{hypergeometric proposition} shows that the
results in this paper easily translate into the language of
hypergeometric functions, the relationship to Mahler measure is more
important than simple pedagogy.  Bertin proved that for certain
values of $u$ the zero varieties of the (projectivized) polynomials
from equations \eqref{bertin g1 definition} and \eqref{bertin g2
definition} define $K3$ hypersurfaces. She also proved formulas
relating the $L$-functions of these $K3$ surfaces at $s=3$ to
rational multiples of the Mahler measures.  Proposition
\ref{hypergeometric proposition} shows that these results imply
explicit $_5F_4$ evaluations (see Corollary
\ref{corollary:hypergeometric explicit evaluations} for explicit
examples). While it might also be interesting to interpret the
polynomials from equations \eqref{f2 definition} through \eqref{f4
definition} in terms of $K3$ hypersurfaces, we will not pursue that
direction here.

\begin{theorem}\label{thm: f2, f3, f4 in terms of G} For $|q|$ sufficiently small
\begin{align}
f_2(s_2(q))=&-\frac{2}{15}G(q)-\frac{1}{15}G\left(-q\right)+\frac{3}{5}G\left(q^2\right),
\label{f2 in terms of G}\\
f_3(s_3(q))=&-\frac{1}{8}G(q)+\frac{3}{8}G\left(q^3\right),\label{f3
in terms of G}\\
f_4(s_4(q))=&-\frac{1}{3}G(q)+\frac{2}{3}G\left(q^2\right),\label{f4
in terms of G}
\end{align}
where
\begin{align*}
s_2(q)=&q^{-1}\left(-q;q^2\right)_{\infty}^{24},\\
s_3(q)=&\frac{1}{q}\left(27q\frac{\left(q^3;q^3\right)_{\infty}^6}{\left(q;q\right)_{\infty}^6}+\frac{\left(q;q\right)_{\infty}^6}{\left(q^3;q^3\right)_{\infty}^6}\right)^2,\\
s_4(q)=&\frac{1}{q}\frac{\left(q^2;q^2\right)_{\infty}^{24}}{\left(q;q\right)_{\infty}^{24}}\left(16q\frac{(q;q)_{\infty}^4\left(q^4;q^4\right)_{\infty}^{8}}{\left(q^2;q^2\right)_{\infty}^{12}}+\frac{\left(q^2;q^2\right)_{\infty}^{12}}{(q;q)_{\infty}^4\left(q^4;q^4\right)_{\infty}^{8}}\right)^4.
\end{align*}
The following inverse relations hold for $|q|$ sufficiently small:
\begin{align}
G(q)=&-19f_2\left(s_2(q)\right)-4f_2\left(s_2(-q)\right)+24f_2\left(s_2\left(q^2\right)\right)-12f_2\left(s_2\left(-q^2\right)\right),\label{G
in terms of f2}\\
\begin{split}
G(q)=&-\frac{19}{2}f_3\left(s_3(q)\right)-\frac{3}{2}f_3\left(s_3\left(\e^{2\pi\ii/3}
q\right)\right)\\&-\frac{3}{2}f_3\left(s_3\left(\e^{4\pi\ii/3}
q\right)\right)+\frac{9}{2}f_3\left(s_3\left(q^3\right)\right),\label{G
in terms of f3}
\end{split}\\
G(q)=&-5f_4\left(s_4(q)\right)-2f_4\left(s_4(-q)\right)+4f_4\left(s_4\left(q^2\right)\right).\label{G
in terms of f4}
\end{align}
\end{theorem}
\begin{proof}We can use Ramanujan's theory of elliptic functions to verify
the first half of this theorem.  Recall that the \textit{elliptic
nome} is defined by
\begin{equation*}
q_j(\alpha)=\exp\left(-\frac{\pi}{\sin\left(\pi/
j\right)}\frac{{_2F_1}\left(\frac{1}{j},1-\frac{1}{j};1;1-\alpha\right)}{{_2F_1}\left(\frac{1}{j},1-\frac{1}{j};1;\alpha\right)}\right).
\end{equation*}
It is a well established fact that
$s_j\left(q_j\left(\alpha\right)\right)$ is a rational function of
$\alpha$ whenever $j\in\{2,3,4\}$. For example if $q=q_2(\alpha)$,
then $s_2\left(q\right)=\frac{16}{\alpha(1-\alpha)}$.  Therefore we
can verify Eq. \eqref{f2 in terms of G} by differentiating with
respect to $\alpha$, and by showing that the identity holds when
$q\rightarrow 0$.

Observe that when $q\rightarrow 0$ both sides of Eq. \eqref{f2 in
terms of G} approach $-\log|q|+O(q)$. Differentiating with respect
to $\alpha$ yields
\begin{equation*}
\begin{split}
-\frac{(1-2\alpha)}{2\alpha(1-\alpha)}{_3F_2}&\left(\substack{\frac{1}{2},\frac{1}{2},\frac{1}{2}\\1,1};4\alpha(1-\alpha)\right)\\
&=-\frac{1}{2q}\left(1-16\sum_{n=1}^{\infty}n^3\frac{
q^n}{1-q^n}+256\sum_{n=1}^{\infty}n^3\frac{
q^{4n}}{1-q^{4n}}\right)\frac{\d q}{\d \alpha}.
\end{split}
\end{equation*}
This final identity follows from applying three well known formulas:
\begin{align*}
&\frac{\d q}{\d
\alpha}=\frac{q}{\alpha(1-\alpha){_2F_1}^2\left(\frac{1}{2},\frac{1}{2};1;\alpha\right)},\\
&{_3F_2}\left(\substack{\frac{1}{2},\frac{1}{2},\frac{1}{2}\\1,1};4\alpha(1-\alpha)\right)={_2F_1}^2\left(\frac{1}{2},\frac{1}{2};1;\alpha\right),\\
&1-16\sum_{n=1}^{\infty}n^3\frac{q^n}{1-q^n}
+256\sum_{n=1}^{\infty}n^3\frac{q^{4n}}{1-q^{4n}}=(1-2\alpha){_2F_1}^4\left(\frac{1}{2},\frac{1}{2};1;\alpha\right).
\end{align*}
We can verify equations \eqref{f3 in terms of G} and \eqref{f4 in
terms of G} in a similar manner by using the fact that
$s_3\left(q_3(\alpha)\right)=\frac{27}{\alpha(1-\alpha)}$, and
$s_4\left(q_4(\alpha)\right)=\frac{64}{\alpha(1-\alpha)}$.

The crucial observation for proving equations \eqref{G in terms of
f2} through \eqref{G in terms of f4} is the fact that $G(q)$
satisfies the following functional equation for any prime $p$:
\begin{equation}\label{G functional equation}
\sum_{j=0}^{p-1}G\left(\e^{2\pi \ii j/p}
q\right)=\left(1+p^3\right)G\left(q^p\right)-p^2
G\left(q^{p^2}\right).
\end{equation}
We will only need the $p=2$ case to prove Eq. \eqref{G in terms of
f2}:
\begin{equation*}\label{G(q) p=2 functional equation}
G(q)+G(-q)=9G\left(q^2\right)-4G\left(q^4\right).
\end{equation*}
Notice that this last formula always allows us to eliminate
$G\left(q^4\right)$ from an equation.  Applying the substitutions
$q\rightarrow-q$, $q\rightarrow q^2$, and $q\rightarrow -q^2$ to Eq.
\eqref{f2 in terms of G} yields
\begin{equation*}
\begin{pmatrix}
&-2/15&\ &-1/15&\ &3/5&\ \\
&-1/15&\ &-2/15&\ &3/5&\ \\
&-3/20&\ &-3/20&\ &23/20&\
\end{pmatrix}
\begin{pmatrix}
G\left(q\right)\\
G\left(-q\right)\\
G\left(q^2\right)
\end{pmatrix}
=
\begin{pmatrix}
f_2\left(s_2(q)\right)\\
f_2\left(s_2(-q)\right)\\
2f_2\left(s_2\left(q^2\right)\right)-f_2\left(s_2\left(-q^2\right)\right)
\end{pmatrix}.
\end{equation*}
Since this system of equations is non-singular, we can invert the
matrix to recover Eq. \eqref{G in terms of f2}.  We can prove
equations \eqref{G in terms of f3} and \eqref{G in terms of f4} in a
similar fashion.$\blacksquare$
\end{proof}

If we compare Theorem \ref{thm: f2, f3, f4 in terms of G} with
Bertin's results we can deduce some obvious relationships between
the Mahler measures.  For example, combining Eq. \eqref{f4 in terms
of G} with Eq. \eqref{bertin g2 series}, and combining Eq. \eqref{f3
in terms of G} with Eq. \eqref{bertin g1 series simplified}, we find
that
\begin{align}
g_1(t_1(q))=&\frac{1}{20}f_4\left(s_4(q)\right)+\frac{3}{20}f_4\left(s_4\left(q^3\right)\right),\label{g1 in terms of f4}\\
g_2\left(t_2(q)\right)=&-\frac{1}{15}f_3\left(s_3(q)\right)
+\frac{8}{15}f_3\left(s_3\left(q^2\right)\right).\label{g2 in terms
of f3}
\end{align}
Notice that many more identities follow from substituting equations
\eqref{G in terms of f2} through \eqref{G in terms of f4} into
formulas \eqref{bertin g1 series simplified}, \eqref{bertin g2
series}, \eqref{f2 in terms of G}, \eqref{f3 in terms of G} and
\eqref{f4 in terms of G}. However, for the remainder of this section
we will restrict our attention to equations \eqref{g1 in terms of
f4} and \eqref{g2 in terms of f3}. In particular, we will appeal to
the theory of elliptic functions to transform those results into
identities which depend on rational arguments.

If we let $q=q_2(\alpha)$, then it is well known that
$\frac{q^{j/24}\left(q^j;q^j\right)_\infty}{q^{1/24}(q;q)_\infty}$
is an algebraic function of $\alpha$ for $j\in \{1,2,3,\dots\}$ (for
example see \cite{Be3} or \cite{Be5}). It follows immediately that
$s_2(q)$, $s_3(q)$, $s_4(q)$, $t_1(q)$, and $t_2(q)$ are also
algebraic functions of $\alpha$. The following lemma lists several
instances where those functions have rational parameterizations.

\begin{lemma}\label{parameterization lemma} Suppose that $q=q_2(\alpha)$, where $\alpha=p(2+p)^3/(1+2p)^3$.
The following identities hold for $|p|$ sufficiently small:
{\allowdisplaybreaks
\begin{align*}
s_2(q)=&\frac{16(1+2p)^6}{p(1-p)^3(1+p)(2+p)^3},&
s_2\left(q^3\right)=&\frac{16(1+2p)^2}{p^3(1-p)(1+p)^3(2+p)},\\
s_2(-q)=&-\frac{16(1-p)^6(1+p)^2}{p(2+p)^3(1+2p)^3},
&s_2\left(-q^3\right)=&-\frac{16(1-p)^2(1+p)^6}{p^3(2+p)(1+2p)},\\
s_2\left(-q^2\right)=&\frac{16^2(1-p)^3(1+p)(1+2p)^3}{p^2(2+p)^6},
&s_2\left(-q^6\right)=&\frac{16^2(1-p)(1+p)^3(1+2p)}{p^6(2+p)^2},\\
s_3(q)=&\frac{4\left(1+4p+p^2\right)^6}{p\left(1-p^2\right)^4(2+p)(1+2p)},&
s_3\left(q^2\right)=&\frac{16\left(1+p+p^2\right)^6}{p^2\left(1-p^2\right)^2(2+p)^2(1+2p)^2},\\
s_3(-q)=&-\frac{4\left(1-2p-2p^2\right)^6}{p(1-p^2)(2+p)(1+2p)^4},&
s_3\left(q^4\right)=&\frac{4\left(2+2p-p^2\right)^6}{p^4\left(1-p^2\right)(2+p)^4(1+2p)},
\end{align*}
\begin{align*}
s_4(q)=&\frac{16\left(1+14p+24p^2+14p^3+p^4\right)^4}{p(1-p)^6(1+p)^2(2+p)^3(1+2p)^3},\\
s_4\left(q^3\right)=&\frac{16\left(1+2p+2p^3+p^4\right)^4}{p^3(1-p)^2(1+p)^6(2+p)(1+2p)},\\
s_4(-q)=&-\frac{16\left(1-10p-12p^2-4p^3-2p^4\right)^4}{p(1-p)^3(1+p)(1+2p)^6(2+p)^3},\\
s_4\left(-q^3\right)=&-\frac{16\left(1+2p-4p^3-2p^4\right)^4}{p^3(1-p)(1+p)^3(1+2p)^2(2+p)}.
\end{align*}
Rational formulas also exist for certain values of $t_1^2(q)$ and
$t_2(q)$:
\begin{align*}
t_1^2(q)=&\frac{4\left(1+p+p^2\right)^2\left(1+4p+p^2\right)^2}{p(1-p^2)^2(2+p)(1+2p)},\\
t_1^2(-q)=&-\frac{4\left(1+p+p^2\right)^2\left(1-2p-2p^2\right)^2}{p(1-p^2)(2+p)(1+2p)^2},\\
t_2(q)=&-\frac{4\left(1-p^2\right)^2}{p(2+p)(1+2p)},\\
t_2(-q)=&-\frac{4(1+p+p^2)^2}{p(1-p^2)(2+p)}.
\end{align*}}
\end{lemma}
The main difficulty with Lemma \ref{parameterization lemma} is the
fact that very few values of $s_j(\pm q^n)$ reduce to rational
functions of $p$. Consider the set $\left\{s_2(q),
s_2(-q),s_2\left(-q^2\right),s_2\left(q^2\right)\right\}$ as an
example. While Lemma \ref{parameterization lemma} shows that
$s_2(q)$, $s_2(-q)$, and $s_2\left(-q^2\right)$ are all rational
with respect to $p$, the formula for $s_2\left(q^2\right)$ involves
radicals. Recall that if $\alpha=p(2+p)^3/(1+2p)^3$, then
\begin{equation*}
s_2\left(q^2\right)=\frac{4\left(1+\sqrt{1-\alpha}\right)^6}{\alpha^2\sqrt{1-\alpha}},
\end{equation*}
where $\sqrt{1-\alpha}=\frac{1-p}{(1+2p)^2}\sqrt{(1-p^2)(1+2p)}.$
Since the curve $X^2=(1-p^2)(1+2p)$ is elliptic with conductor $24$,
it follows immediately that rational substitutions for $p$ will
never reduce $s_2\left(q^2\right)$ to a rational function.  For the
sake of legibility, we will therefore avoid all identities which
involve those four functions simultaneously.  By avoiding pitfalls
of this nature, we can derive several interesting results from Lemma
\ref{parameterization lemma}.

\begin{theorem}\label{theorem: g1 and g2 rational transformations} For $|z|$ sufficiently large
\begin{align}
g_1\left(3\left(z+z^{-1}\right)\right)=&\frac{1}{20}f_4\left(\frac{9\left(3+z^2\right)^4}{z^6}\right)+\frac{3}{20}f_4\left(\frac{9\left(3+z^{-2}\right)^4}{z^{-6}}\right)\label{g1 in terms of f4 explicit},\\
g_2(z)=&-\frac{1}{15}f_3\left(\frac{(16-z)^3}{z^2}\right)+\frac{8}{15}f_3\left(-\frac{(4-z)^3}{z}\right).\label{g2
in terms of f3 explicit}
\end{align}
\end{theorem}
\begin{proof} These identities follow
from applying Lemma \ref{parameterization lemma} to equations
\eqref{g1 in terms of f4} and \eqref{g2 in terms of f3}.  If we
consider Eq. \eqref{g1 in terms of f4}, then Lemma
\ref{parameterization lemma} shows that $t_1^2(q)$, $s_4(q)$, and
$s_4\left(q^3\right)$ are all rational functions of $p$. Forming a
resultant with respect to $p$, we obtain
\begin{equation*}
\begin{split}
0=\mathop{\Res}_p&\left[\frac{4\left(1+p+p^2\right)^2\left(1+4p+p^2\right)^2}{p\left(1-p^2\right)^2(2+p)(1+2p)}-t_1^2(q),\right.\\
&\quad\left.\frac{16\left(1+14p+24p^2+14p^3+p^4\right)^4}{p(1-p)^6(1+p)^2(2+p)^3(1+2p)^3}-s_4(q)\right].
\end{split}
\end{equation*}
Simplifying with the aid of a computer, this becomes
\begin{equation*}
0=s_4^2(q)+\left(12+t_1^2(q)\right)^4-s_4(q)\left(-288+352t_1^2(q)-42t_1^4(q)+t_1^6(q)\right).
\end{equation*}
If we choose $z$ so that $t_1(q)=3\left(z+z^{-1}\right)$, then
$s_4(q)=9\left(3+z^2\right)^4 z^{-6}$, and a formula for
$s_4\left(q^3\right)$ follows in a similar fashion.$\blacksquare$
\end{proof}
%
If we let $u=1/z$ with $z\in\mathbb{R}$ and sufficiently large, then
Eq. \eqref{g2 in terms of f3 explicit} reduces to the following
infinite series identity:
\begin{equation}\label{5F4 identity one}
\begin{split}
 \sum_{n=1}^{\infty}\frac{u^{n}}{n}&\sum_{k=0}^{n}{2k\choose
k}{2n-2k\choose n-k}{n\choose
k}^2\\=&\frac{1}{5}\log\left(\frac{(1-16u)}{(1-4u)^8}\right)
+\frac{4u}{5(1-16u)^3}{_5F_4}\left(\substack{\frac{4}{3},\frac{3}{2},\frac{5}{3},1,1\\2,2,2,2};-\frac{108u}{(1-16u)^3}\right)\\
&+\frac{32u^2}{5(1-4u)^3}{_5F_4}\left(\substack{\frac{4}{3},\frac{3}{2},\frac{5}{3},1,1\\2,2,2,2};\frac{108u^2}{(1-4u)^3}\right).
\end{split}
\end{equation}
Similarly, if we let $u=1/z^2$ then Eq. \eqref{g1 in terms of f4
explicit} is equivalent to:
\begin{equation}\label{5F4 identity two}
\begin{split}
\sum_{n=1}^{\infty}\frac{1}{n}&\left(\frac{u}{9(1+u)^2}\right)^n{2n\choose
n}\sum_{k=0}^{n}{2k\choose k}{n\choose k}^2\\
=&\frac{2}{5}\log\left(\frac{27(1+u)^5}{(3+u)^3(1+3u)}\right)
+\frac{4u^3}{5(3+u)^4}{_5F_4}\left(\substack{\frac{5}{4},\frac{3}{2},\frac{7}{4},1,1\\2,2,2,2};\frac{256u^3}{9(3+u)^4}\right)\\
&+\frac{4u}{15(1+3u)^4}{_5F_4}\left(\substack{\frac{5}{4},\frac{3}{2},\frac{7}{4},1,1\\2,2,2,2};\frac{256u}{9(1+3u)^4}\right).
\end{split}
\end{equation}
In Section \ref{pi section} we will differentiate equations
\eqref{5F4 identity one} and \eqref{5F4 identity two} to obtain
several new formulas for $1/\pi$.  But first we will conclude this
section by deducing some explicit $_5F_4$ evaluations.

    Recall that for certain values of $u$, Bertin evaluated $g_1(u)$ and $g_2(u)$
in terms of the $L$-series of $K3$ surfaces.  She also proved
equivalent formulas involving twisted cusp forms. Amazingly, her
formulas correspond to cases where the right-hand sides of equations
\eqref{g1 in terms of f4 explicit} and \eqref{g2 in terms of f3
explicit} collapse to one hypergeometric term. We can combine her
results with equations \eqref{5F4 identity one} and \eqref{5F4
identity two} to deduce several new $_5F_4$ evaluations.
\begin{corollary}\label{corollary:hypergeometric explicit evaluations}If $g(q)=q\left(q^2;q^2\right)_{\infty}^3\left(q^6;q^6\right)_\infty^3$,
then
\begin{equation}\label{formula 1 for 5F4 cusp form at s=3}
{_5F_4}\left(\substack{\frac{4}{3},\frac{3}{2},\frac{5}{3},1,1\\2,2,2,2};1\right)=18\log(2)+27\log(3)-\frac{810\sqrt{3}}{\pi^3}L(g,3).
\end{equation}
If
$f(q)=q\left(q;q\right)_{\infty}^2\left(q^2;q^2\right)_{\infty}\left(q^4;q^4\right)_{\infty}\left(q^8;q^8\right)_{\infty}^2$,
then
\begin{equation}\label{formula 2 for 5F4 cusp form at s=3}
{_5F_4}\left(\substack{\frac{5}{4},\frac{3}{2},\frac{7}{4},1,1\\2,2,2,2};1\right)=\frac{256}{3}\log(2)-\frac{5120\sqrt{2}}{3\pi^3}L(f,3).
\end{equation}
\end{corollary}

While many famous $_5F_4$ identities, such as Dougall's formula
\cite{Be2}, reduce special values of the $_5F_4$ function to gamma
functions, equations \eqref{formula 1 for 5F4 cusp form at s=3} and
\eqref{formula 2 for 5F4 cusp form at s=3} do not fit into this
category.  Rather these new formulas are higher dimensional
analogues of Boyd's conjectures.  In particular, Boyd has
conjectured large numbers of identities relating two-variable Mahler
measures (that mostly reduce to $_4F_3$ functions) to the $L$-series
of elliptic curves \cite{Bo1}. The most famous outstanding
conjecture of this type asserts that
\begin{equation*}
\m\left(1+x+\frac{1}{x}+y+\frac{1}{y}\right)=-2\Re\left[{_4F_3}\left(\substack{\frac{3}{2},\frac{3}{2},1,1\\2,2,2};16\right)\right]\stackrel{?}{=}\frac{15}{4\pi^2}L(f,2),
\end{equation*}
where
\begin{equation*}
f(q)=q\prod_{n=1}^{\infty}\left(1-q^n\right)\left(1-q^{3n}\right)\left(1-q^{5n}\right)\left(1-q^{15n}\right),
\end{equation*}
and ``$\displaystyle\stackrel{?}{=}$" indicates numerical equality
to at least $50$ decimal places.  Recently, Kurokawa and Ochiai
proved a formula \cite{KO} which simplifies this last conjecture to
\begin{equation*}
{_3F_2}\left(\substack{\frac{1}{2},\frac{1}{2},\frac{1}{2}\\\frac{3}{2},1};\frac{1}{16}\right)\stackrel{?}{=}\frac{15}{\pi^2}L(f,2).
\end{equation*}
Of course it would be highly desirable to rigorously prove Boyd's
conjectures. Failing that, it might be interesting to search for
more hypergeometric identities like equations \eqref{formula 1 for
5F4 cusp form at s=3} and \eqref{formula 2 for 5F4 cusp form at
s=3}.  This line of thought suggests the following fundamental
problem with which we shall conclude this section:

\bigskip
\noindent\textbf{Open Problem:} Determine every $L$-series that can
be expressed in terms of generalized hypergeometric functions with
algebraic parameters.
\bigskip
\section{New formulas for $1/\pi$}
\label{pi section} \init

In the previous section we produced several new transformations for
the $_5F_4$ hypergeometric function.  Now we will differentiate
those formulas to obtain some new $_3F_2$ transformations, and
several accompanying formulas for $1/\pi$. The following formula is
a typical example of the identities in this section:
\begin{equation}\label{example 1/pi formula}
\frac{2}{\pi}=\sum_{n=0}^{\infty}(-1)^n\frac{(3n+1)}{32^n}\sum_{k=0}^{n}{2n-2k\choose
n-k}{2k\choose k}{n\choose k}^2.
\end{equation}
Ramanujan first proved identities like Eq. \eqref{example 1/pi
formula} in his famous paper ``Modular equations and approximations
to $\pi$" \cite{Ra}. He showed that the following infinite series
holds for certain constants $A$, $B$, and $X$:
\begin{equation}\label{ramanujan 1/pi formula}
\frac{1}{\pi}=\sum_{n=0}^{\infty}(A
n+B)\frac{\left(1/2\right)^3_n}{n!^3}X^n
\end{equation}
Ramanujan determined many sets of algebraic values for $A$, $B$, and
$X$ by expressing them in terms of the classical singular moduli
$G_n$ and $g_n$. He also stated (but did not prove) several formulas
for $1/\pi$ where $\left(1/2\right)^3_n$ is replaced by
$\left(1/a\right)_n\left(1/2\right)_n\left(1-1/a\right)_n$ for
$a\in\left\{3,4,6\right\}$ (for more details see \cite{BR} or
\cite{Ch}).

Ramanujan's formulas for $1/\pi$ have attracted a great deal of
attention because of their intrinsic beauty, and because they
converge extremely quickly.  For example, \textit{Mathematica}
calculates $\pi$ using a variant of a Ramanujan-type formula due to
the
Chudnovsky brothers \cite{mw}: 
\begin{equation}\label{chudnovsky example}
\frac{1}{\pi}=12\sum_{n=0}^{\infty}\frac{(-1)^n(6n)!(13591409+54513013n)}{n!^3(3n)!\left(640320^3\right)^{n+1/2}}.
\end{equation}
More recent mathematicians including Yang and Zudilin have derived
formulas for $1/\pi$ which are not hypergeometric, but still similar
to Eq. \eqref{chudnovsky example}. For example, Yang showed that
\begin{equation*}
\frac{18}{\pi\sqrt{15}}=\sum_{n=0}^{\infty}\frac{(4n+1)}{36^n}\sum_{k=0}^{n}{n\choose
k}^4,
\end{equation*}
and Zudilin gave many infinite series for $1/\pi$ containing nested
sums of binomial coefficients \cite{Zd}.  All of the formulas that
we will prove, including Eq. \eqref{example 1/pi formula}, are
essentially of this type.  Before proving the next theorem, we will
point out that equation \eqref{transformation 3f2 number 1} appears
implicitly in the work of Chan, Chan and Liu, and can be obtained by
combining equations (4.5) and (4.6) in their paper \cite{ccl}.

\begin{theorem}\label{thm: 3f2 transformations} For $|u|$ sufficiently small
\begin{equation}\label{transformation 3f2 number 1}
{_3F_2}\left(\substack{\frac{1}{3},\frac{1}{2},\frac{2}{3}\\1,1};\frac{108u^2}{(1-4u)^3}\right)=(1-4u)\sum_{n=0}^{\infty}u^n
\sum_{k=0}^{n}{2n-2k \choose n-k}{2k\choose k}{n\choose k}^2.
\end{equation}
If $|u|$ is sufficiently small
\begin{equation}\label{transformation 3f2 number 2}
{_3F_2}\left(\substack{\frac{1}{4},\frac{1}{2},\frac{3}{4}\\1,1};\frac{256u}{9(1+3u)^4}\right)=\frac{(1+3u)}{(1+u)}\sum_{n=0}^{\infty}\left(\frac{u}{9(1+u)^2}\right)^n{2n\choose
n}\sum_{k=0}^{n}{2k\choose k}{n\choose k}^2.
\end{equation}
\end{theorem}
\begin{proof} Applying the operator $u\frac{\d}{\d u}$ to Eq.
\eqref{5F4 identity one}, and then simplifying yields
\begin{equation*}
\begin{split}
\sum_{n=0}^{\infty}u^n\sum_{k=0}^{n}{2k\choose k}{2n-2k\choose
n-k}{n\choose
k}^2=&-\frac{(1+32u)}{15(1-16u)}{_3F_2}\left(\substack{\frac{1}{3},\frac{1}{2},\frac{2}{3}\\1,1};-\frac{108u}{(1-16u)^3}\right)\\
&+\frac{16(1+2u)}{15(1-4u)}{_3F_2}\left(\substack{\frac{1}{3},\frac{1}{2},\frac{2}{3}\\1,1};\frac{108u^2}{(1-4u)^3}\right).
\end{split}
\end{equation*}
Eq. \eqref{transformation 3f2 number 1} then follows from applying a
$_3F_2$ transformation:
\begin{equation}\label{aux 3f2 trans 1}
{_3F_2}\left(\substack{\frac{1}{3},\frac{1}{2},\frac{2}{3}\\1,1};-\frac{108u}{(1-16u)^3}\right)
=\frac{(1-16u)}{(1-4u)}{_3F_2}\left(\substack{\frac{1}{3},\frac{1}{2},\frac{2}{3}\\1,1};\frac{108u^2}{(1-4u)^3}\right).
\end{equation}
We can prove Eq. \eqref{transformation 3f2 number 2} in a similar
manner by differentiating Eq. \eqref{5F4 identity two} and then
using
\begin{align}
{_3F_2}\left(\substack{\frac{1}{4},\frac{1}{2},\frac{3}{4}\\1,1};\frac{256u^3}{9(3+u)^4}\right)
=&\frac{(3+u)}{3(1+3u)}{_3F_2}\left(\substack{\frac{1}{4},\frac{1}{2},\frac{3}{4}\\1,1};\frac{256u}{9(1+3u)^4}\right).\label{aux
3f2 trans 2}
\end{align}
Equation \eqref{aux 3f2 trans 1} can each be derived in three steps.
First square both sides of the following $_2F_1$ identity (see
Corollary 6.2 in \cite{Be5}):
\begin{equation*}
{_2F_1}\left(\substack{\frac{1}{3}\frac{2}{3}\\1};\frac{(1-p)(2+p)^2}{4}\right)
=\frac{2}{(1+p)}{_2F_1}\left(\substack{\frac{1}{3},\frac{2}{3}\\1};\frac{(1-p)^2(2+p)}{2(1+p)^3}\right).
\end{equation*}
When $|1-p|$ is sufficiently small, we can apply the following
$_3F_2$ transformation
\begin{equation*}
{_3F_2}\left(\substack{\frac{1}{3},\frac{1}{2},\frac{2}{3}\\
1,1};4x(1-x)\right)={_2F_1}^2\left(\substack{\frac{1}{3},\frac{2}{3}\\
1}; x\right),
\end{equation*}
and then conclude by setting $u=\frac{(-1+p)(3+p)}{16p(2+p)}$.
Similarly, equation \eqref{aux 3f2 trans 2} follows from combining
another $_2F_1$ identity (which can be shown to follow from Theorem
9.15 in \cite{Be5}):
\begin{equation*}
{_2F_1}\left(\substack{\frac{1}{4},\frac{3}{4}\\ 1};
1-\frac{64p}{(3+6p-p^2)^2}\right)=\sqrt{\frac{3+6p-p^2}{27-18p-p^2}}{_2F_1}\left(\substack{\frac{1}{4},\frac{3}{4}\\
1};1-\frac{64p^3}{(27-18p-p^2)^2}\right),
\end{equation*}
with a similar $_3F_2$ transformation
\begin{equation*}
{_3F_2}\left(\substack{\frac{1}{4},\frac{1}{2},\frac{3}{4}\\
1,1};4x(1-x)\right)={_2F_1}^2\left(\substack{\frac{1}{4},\frac{3}{4}\\
1}; x\right),
\end{equation*}
and then setting $u=\frac{9(1-p)}{p(9-p)}$. $\blacksquare$
\end{proof}
    While the infinite series in Theorem \ref{thm:
3f2 transformations} are not hypergeometric since they involve
nested binomial sums, they are still interesting.  In particular,
those formulas easily translate into unexpected integrals involving
powers of modified Bessel functions.  For $|x|$ sufficiently small,
Eq. \eqref{transformation 3f2 number 1} is equivalent to
\begin{equation}\label{I0^3 integral}
\int_{0}^{\infty}\e^{-3\left(x+x^{-1}\right)u}I_0^3\left(2u\right)\d
u=\frac{x}{3\left(1+3x^2\right)}{_3F_2}\left(\substack{\frac{1}{4},\frac{1}{2},\frac{3}{4}\\1,1}\frac{256x^2}{9(1+3x^2)^4}\right),
\end{equation}
where $I_0(u)$ is the modified Bessel function of the first kind.
Recall the series expansions for $I_0(2u)$ and $I_0^2(2u)$:
\begin{align*}
I_0(2u)=\sum_{n=0}^{\infty}\frac{u^{2n}}{n!^2},&&I_0^2(2u)=\sum_{n=0}^{\infty}{2n\choose
n}\frac{u^{2n}}{n!^2}.
\end{align*}
Eq. \eqref{I0^3 integral} is surprising because there is no known
hypergeometric expression $I_0^3(2u)$ \cite{bender}.  It is
therefore not obvious that the Laplace transform of $I_0^3(2u)$
should equal a hypergeometric function.  M. Lawrence Glasser has
kindly pointed out that equation \eqref{I0^3 integral} is
essentially a well known result, and that a variety of similar
integrals have also been studied by Joyce \cite{Jc}, Glasser and
Montaldi \cite{Gl}, and others.

Finally, we will list a few formulas for $1/\pi$. Notice that
equation \eqref{pi 2} first appeared in the work of Chan, Chan and
Liu \cite{ccl}.

\begin{corollary}\label{pi cor} Let $a_n=\sum_{k=0}^{n}{2n-2k\choose
n-k}{2k\choose k}{n\choose k}^2$, then the following formulas are
true:
\begin{align}
\frac{2}{\pi}=&\sum_{n=0}^{\infty}(-1)^n\frac{(3n+1)}{32^{n}}a_n,\label{pi 1}\\
\frac{8\sqrt{3}}{3\pi}=&\sum_{n=0}^{\infty}\frac{(5n+1)}{64^{n}}a_n,\label{pi 2}\\
\frac{9+5\sqrt{3}}{\pi}=&\sum_{n=0}^{\infty}(6n+3-\sqrt{3})\left(\frac{3\sqrt{3}-5}{4}\right)^n
a_n.\label{pi 3}
\end{align}
Let $b_n={2n\choose n}\sum_{k=0}^{n}{2k\choose k}{n\choose k}^2$,
then the following identity holds:
\begin{align}\label{pi 4}
\frac{2\left(64+29\sqrt{3}\right)}{\pi}=&\sum_{n=0}^{\infty}\left(520n+159-48\sqrt{3}\right)\left(\frac{80\sqrt{3}-139}{484}\right)^n
b_n.
\end{align}
\end{corollary}
\begin{proof} We can use Eq. \eqref{transformation 3f2 number 1} to easily deduce that if
\begin{equation*}
\begin{split}
\sum_{n=0}^{\infty}&(a
n+b)\frac{(1/3)_n(1/2)_n(2/3)_n}{n!^3}\left(\frac{108u^2}{(1-4u)^3}\right)^n
=\sum_{n=0}^{\infty}(A n+B)u^n a_n,
\end{split}
\end{equation*}
then $A=a(1-4u)/(2+4u)$, and $B=a(-4u)(1-4u)/(2+4u)+b(1-4u)$. Since
the left-hand side of this last formula equals $1/\pi$ when
$\left(a,b,\frac{108u^2}{(1-4u)^3}\right)\in\{\left(\frac{60}{27},\frac{8}{27},\frac{2}{27}\right),\left(\frac{2}{\sqrt{3}},\frac{1}{3\sqrt{3}},\frac{1}{2}\right),(\frac{45}{11}-\frac{5}{33}\sqrt{3},\frac{6}{11}-\frac{13}{99}\sqrt{3},-\frac{194}{1331}+\frac{225}{2662}\sqrt{3})\}$,
it is easy to verify equations \eqref{pi 1} through \eqref{pi 3}
\cite{Ch}.

    We can verify Eq. \eqref{pi 4} in a similar manner by combining Eq. \eqref{transformation 3f2 number 2} with Ramanujan's formula
\begin{equation*}
\frac{8}{\pi}=\sum_{n=0}^{\infty}(20
n+3)\frac{(1/4)_n(1/2)_n(3/4)_n}{n!^3}\left(\frac{-1}{4}\right)^n.
\end{equation*}$\blacksquare$
\end{proof}
\section{Conclusion}
\label{conclusion}\init
We will conclude the paper by suggesting two future projects.
Firstly, it would be desirable to determine whether or not a
rational series involving $b_n$ exists for $1/\pi$.  Secondly, it
might be interesting to consider the Mahler measure
\begin{equation*}
f_6(u)=\m\left(u-\left(z+z^{-1}\right)^6\left(y+y^{-1}\right)^2(1+x)^3x^{-2}\right),
\end{equation*}
since $f_6(u)$ arises from Ramanujan's theory of signature $6$.

\bigskip
\begin{acknowledgements}

The author would like to thank David Boyd for the many helpful
discussions and encouragement. The author also thanks Fernando
Rodriquez-Villegas and Marie Jos\'{e} Bertin for the useful
discussions.  The author thanks Wadim Zudilin for the useful
communications and for the reference \cite{Zd}.  The author also
thanks Heng Huat Chan and Zhiguo Liu for the reference and their
useful remarks \cite{ccl}. The author extends his gratitude to M.
Lawrence Glasser for the references \cite{Gl} and \cite{Jc}.
Finally, thanks to Jianyun Shan for pointing out \cite{JZ}.
\end{acknowledgements}


\end{document}